\documentclass{article}
\usepackage{amssymb}
\title{Dynamics of a Compact Operator }
\author{Teck-Cheong Lim \\ Department of Mathematical Sciences \\ George Mason University \\ 4400, University Drive\\Fairfax, VA 22030\\U.S.A.\\
{\it e-mail address}: tlim@gmu.edu}
\newtheorem {co}{Corollary}
\newtheorem {ex}{Example}

\newtheorem {df} {Definition}
\newtheorem {theo} {Theorem}

\newtheorem {rem} {Remark}

\begin{document}
\maketitle

\begin{abstract}
Let $T:X\to X$ be a compact linear (or more generally affine) operator from a Banach space into itself. For each $x\in X$, the sequence of iterates $T^nx, n=0,1,\cdots$ and its averages $\frac{1}{k}\sum_{k=0}^nT^{k-1}x, n=0,1,\cdots$ are either bounded or approach infinity.
\end{abstract}


Let $X$ be a set and $f: X\to X$ a map from $X$ into $X$. For an $x\in X$, the sequence of iterations $x, f(x), f^2(x),\cdots, f^k(x),\cdots$ can be considered as a trajectory of a dynamical system where time is the discrete nonnegative integers: starting with the initial (time $t=0$) state x, the state at time $t=k$ is $f^k(x)$.
Suppose now $X=\mathbb{C}^n$ and $f$ is the transformation defined by an $n\times n$ complex matrix $A$. What can one say about the general behavior of trajectories of A? More generally, we consider affine maps on $X$, i.e. maps of the form $Ax+c$, where $A$ is linear and $c$ is a constant vector, and we also allow $X$ to be infinite dimensional. Moreover, we study the behavior of the sequence of averages:
\[\mbox{Ave}_kf(x)=\frac{1}{k}(x+f(x)+\cdots+f^{k-1}(x)), k=1,2,\cdots\]
Note that the method of averaging was used in [2] in approximating solutions of a system of linear equations. In case of linear operators, problems about the linear span of the iterates $T^kx, n=0,1,\cdots$ can be found in [4].\\
 Recall that a square matrix $N$ is called {\em nilpotent} if $N^s=0$ for some nonnegative integer $s$.\\
 Throughout this paper, for nonnegative integers $k,j, k\ge j$, $C(k,j)$ denotes the binomial coefficient \[\frac{k!}{j!(k-j)!}\]
 By convention, $C(k,0)=1$ for $k\ge 0$ and $C(k,j)=0$ if $k<j$.\\ If $k<j$, the sum $\sum_{i=j}^k u_i$ is considered as an empty sum and its value is 0.
\begin{theo} Let $T: \mathbb{C}^n \to \mathbb{C}^n$ be an affine map defined by $Tx=Ax+c$, where $A$ is an $n\times n$ complex matrix and $c$ a constant vector in $\mathbb{C}^n$. Let $\|\cdot\|$ be a norm on $\mathbb{C}^n$. Then for any vector $x\in\mathbb{C}^n$, the sequence \[T^kx,k=0,1,2,\cdots\] is either bounded or $\lim_k \|T^k x\|=\infty$.
\end{theo}
{\bf Proof.}\\
By Jordan canonical decomposition theorem, $\mathbb{C}^n=V_1\oplus V_2\oplus\cdots\oplus V_m$ for some subspaces $V_i, i=1,2,\cdots,m$ with the following properties: (a) each $V_i$ is an invariant subspace of $A$, i.e. $Av\in V_i$ for all $v\in V_i$, and (b) there exists $\lambda_i\in\mathbb{C}$ and a nilpotent matrix $N_i$ such that $Av=\lambda_iv+N_iv$ for all $v\in V_i$. \\Let $P_i$ be the algebraic projection of $\mathbb{C}^n$ onto $V_i$ associated with the decomposition $\mathbb{C}^n=V_1\oplus V_2\oplus\cdots\oplus V_m$.\\
Define a new norm $|\cdot|$ on $\mathbb{C}^n$ by \[|v|=\|P_1v\|+\|P_2v\|+\cdots+\|P_mv\|.\] \\
Let $x$ be a vector in $\mathbb{C}^n$. Let  \[x_k=T^kx=A^kx+c+Ac+\cdots+A^{k-1}c, k=1,2,\cdots.\]
For any vector $x$, we have the following equalities:
\begin{eqnarray*} x &=& P_1x+\cdots+P_mx\\
A^kx &=& A^kP_1x+\cdots+A^kP_mx\\
A^kx &=& P_1A^kx+\cdots+P_mA^kx
\end{eqnarray*}
The last equality follows from the commutative property of $A$ and $P_i,i=1,\cdots,m$ since $A$ is invariant in each $V_j$.\\ Fix an $i$ and write $v=P_ix, d=P_ic, \lambda=\lambda_i, V=V_i, N=N_i$. Let $s$ be the smallest nonnegative integer such that $N^sv=0$ and $t$ the smallest nonnegative integer such that $N^td=0$. Then for $k>\max\{s,t\}$ and $s,t\ge 1$ (if $s=0$ or $t=0$, the corresponding sum below is defined as 0), one has
\begin{eqnarray}
P_ix_k &=& (\lambda I+N)^kv+d+(\lambda I+N)d+\cdots +(\lambda I+N)^{k-1}d\\
&=& \sum_{j=0}^{s-1}C(k,j)\lambda^{k-j}N^jv +\sum_{j=0}^{t-1}S(j,k)N^jd
\end{eqnarray}
where \[S(j,k)=C(j,j)+C(j+1,j)\lambda+\cdots +C(k-1,j)\lambda^{k-1-j}.\] Using the identity $C(j,i+1)+C(j,i)=C(j+1,i+1)$, we have for $\lambda= 1$, $S(j,k)=C(k,j+1)$, and for $\lambda\neq 1$, $S(j,k)$ is given recursively by \[S(0,k)=\frac{1-\lambda^k}{1-\lambda}\] and (by subtracting $\lambda S(j,k)$ from $S(j,k)$)\[(1-\lambda)S(j,k)=S(j-1,k)-\lambda^{k-j}C(k,j), j=1,2,\cdots,t-1, (\mbox{ for } t\ge 2),\]
from which we get an alternate formula for $S(j,k)$: \begin{equation}S(j,k)=\frac{1-\lambda^k}{(1-\lambda)^{j+1}}-\sum_{i=0}^{j-1}\frac{C(k,i+1)\lambda^{k-i-1}}{(1-\lambda)^{j-i}}, j=1,\cdots,t-1\end{equation}
Note that (3) is also valid for $j=0$ for $\sum_{i=0}^{-1} $ is an empty sum. We shall show that for $\lambda\neq 1$,
\begin{equation}P_ix_k = B+\sum_{j=0}^{w-1} \lambda^{k-j}C(k,j)A_j\end{equation}
where $w=\max\{s,t\}$, and $A_j, B, j=0,\cdots,w-1$ are constant vectors independent of $k$; and for $\lambda=1$,
\begin{equation}P_ix_k = v+\sum_{j=1}^l C(k,j)B_j\end{equation}
where $l=\max\{s-1,t\}$, and $B_j, j=1,\cdots,l$ are constant vectors independent of $k$. More precisely, for $k>\max\{s,t\}$, \[A_j=N^jv-\sum_{i=j}^{t-1}\frac{1}{(1-\lambda)^{i-j+1}}N^{i}d, j=0,\cdots,\min\{s,t\}-1,\] \[A_j=\epsilon(s-t)N^jv-\epsilon(t-s)\sum_{i=j}^{t-1}\frac{1}{(1-\lambda)^{i-j+1}}N^{i}d, j=\min\{s,t\},\cdots,w-1\]and \[B=\sum_{i=0}^{t-1}\frac{1}{(1-\lambda)^{i+1}}N^{i}d;\]
and   \[B_j=N^jv+N^{j-1}d, j=1,\cdots, \min\{s-1,t\},\]\[B_j=\epsilon(s-1-t)N^jv+\epsilon(t-s+1)N^{j-1}d, j= \min\{s-1,t\}+1,\cdots,l,\] where $\epsilon(r)=1$ for $r\ge 0$ and $\epsilon(r)=0$ for $r<0$.\\
Indeed substituting (3) into (2), we get
\begin{eqnarray*}P_ix_k &=&\sum_{j=0}^{s-1}C(k,j)\lambda^{k-j}N^jv+\sum_{j=0}^{t-1}\left(\frac{1-\lambda^k}{(1-\lambda)^{j+1}}-\sum_{i=0}^{j-1}\frac{C(k,i+1)\lambda^{k-i-1}}{(1-\lambda)^{j-i}}\right)N^jd\\
&=&\sum_{j=0}^{s-1}C(k,j)\lambda^{k-j}N^jv+\sum_{j=0}^{t-1}\frac{1}{(1-\lambda)^{j+1}}N^{j}d-\sum_{j=0}^{t-1}\sum_{i=0}^j\frac{C(k,i)\lambda^{k-i}}{(1-\lambda)^{j-i+1}}N^jd\\
&=&\sum_{j=0}^{s-1}C(k,j)\lambda^{k-j}N^jv+B-\sum_{i=0}^{t-1}\sum_{j=i}^{t-1}\frac{C(k,i)\lambda^{k-i}}{(1-\lambda)^{j-i+1}}N^jd\\
&=&\sum_{j=0}^{s-1}C(k,j)\lambda^{k-j}N^jv+B-\sum_{j=0}^{t-1}\sum_{i=j}^{t-1}\frac{C(k,j)\lambda^{k-j}}{(1-\lambda)^{i-j+1}}N^id\\
&=&\sum_{j=0}^{s-1}C(k,j)\lambda^{k-j}N^jv+B-\sum_{j=0}^{t-1}C(k,j)\lambda^{k-j}\sum_{i=j}^{t-1}\frac{1}{(1-\lambda)^{i-j+1}}N^id
\end{eqnarray*}
Now (4) follows by considering cases $s> t, s=t$ or $s< t$.\\
Next consider the case $\lambda=1$. Substituting $S(j,k)=C(k,j+1)$ into (2), we get
\begin{eqnarray*}P_ix_k
&=& \sum_{j=0}^{s-1}C(k,j)N^jv+\sum_{j=0}^{t-1}C(k,j+1)N^jd\\
&=& v+\sum_{j=1}^{s-1}C(k,j)N^jv+\sum_{j=1}^{t}C(k,j)N^{j-1}d
\end{eqnarray*}
Now (5) follows by considering cases $s-1>t, s-1=t$ or $s-1>t$.\\ We now proceed to finish the proof of the theorem:\\
Case 1: $|\lambda|>1$. If one of these $A_j, j=0,\cdots,w-1$ is nonzero, then we see from (4) that $\|P_ix_k\|\to \infty$ as $k\to\infty$; otherwise $P_ix_k=B$ is a constant vector.\\ \\
Case 2: $|\lambda|=1$ and $\lambda\neq 1$. If one of these $A_j, j=1,\cdots,w-1$ is nonzero, then we see from (4) that $\|P_ix_k\|\to \infty$ as $k\to\infty$; otherwise the sequence $P_ix_k=B+\lambda^kA_0, k=1,\cdots,$ is bounded.\\ \\
Case 3: $\lambda=1$. If one of these $B_j, j=1,\cdots,l$ is nonzero, then we see from (5) that $\|P_ix_k\|\to \infty$ as $k\to\infty$; otherwise $P_ix_k$ is the constant vector $v$.\\ \\
Case 4: $|\lambda|<1$. We see from (4) that the sequence $P_ix_k, k=1,\cdots,$ converges to the constant vector $B$.\\ \\
We conclude that for each $i$, the sequence $P_ix_k, k=1,\cdots$ is either bounded or tends to infinity. If one of these sequences tends to infinity, then since $|x_k|\ge \|P_ix_k\|$, the sequence $x_k, k=1,\cdots$ tends to infinity in norm $|\cdot|$ and hence in norm $\|\cdot\|$, since the two norms are equivalent. Otherwise all sequences $P_ix_k, k=1,\cdots$ are bounded for all $i=1,\cdots,m$, and from which it follows that $x_k, k=1,\cdots$ is bounded.
This completes the proof. Q.E.D.\\ \\

The following example shows that in infinite dimensional spaces, Theorem 1 is false.

\begin{ex}\rm For nonnegative integers $n$, let $c_n=\frac{1}{2}n(n+1)$. Define $\lambda_i=\frac{1}{2}$ for
$c_{2n}\le i\le c_{2n+1}-1, n=0,1,\cdots$ and  $\lambda_i=2$ for
$c_{2n-1}\le i\le c_{2n}-1, n=1,2,\cdots$. Define a linear
operator $A: l_2\to l_2$ such that $Ae_i=\lambda_i e_{i+1},
i=0,1,\cdots$. Then the sequence of iterates $A^ke_0=\lambda_0\lambda_1\cdots\lambda_{k-1}e_k,
k=1,2,\cdots$ contains subsequences
that converge to 0, (e.g. $A^{c_{2n-1}}e_0=1/2^n, n=1,2,\cdots$), subsequences
that approach infinity (e.g.$A^{c_{2n}}e_0=2^n, n=1,2,\cdots$), and infinitely many bounded
nonconvergent subsequences (e.g. $A^ne_0=e_n$ for $n=2,4,8,12,18,24,32,40,\cdots$). Clearly $A$ is bounded but not compact.\\
\end{ex}
\begin{rem} \rm For $T$ linear, Theorem 1 appeared in a 1997 unpublished article ``On the behavior of the iterates of a matrix" of the author.
\end{rem}
In the case that the mapping $T$ in Theorem 1 is linear, we can actually say more:\\
\begin{theo} Let $A: \mathbb{C}^n \to \mathbb{C}^n$ be a linear map. Let $\|\cdot\|$ be a norm on $\mathbb{C}^n$. Then for any vector $x\in\mathbb{C}^n$, either $\lim_{k\to\infty}T^kx=0$ or  $\lim_{k\to\infty} \|T^k x\|=\infty$ or $F\le\|T^kx\|\le G$ for sufficiently large $k$'s, where $F,G$ are positive numbers with $ F\le G$.
\end{theo}
{\bf Proof.}\\
If $c=0$ in Theorem 1, then the constant vectors $B$ and $d$ in its proof are 0. It follows that  in all cases where $P_ix_k, k=1,2,\ldots$ is bounded and not convergent to 0, it is bounded away from 0. If $P_ix_k$ is bounded away from 0 for some $i$, then so is $x_k=T^kx$ since $|x_k|\ge \|P_ix_k\|$. Otherwise $P_ix_k\to 0$ converges to 0 for all $i$ and hence so does $x_k$. Q.E.D.\\ \\
\begin{ex}\rm The following simple example shows that linearity is needed in Theorem 2: Let $T:\mathbb{C}\to\mathbb{C}$ be the map $Tx=ix+c$, where $c\in\mathbb{C}$ is nonzero and $i=\sqrt{-1}$. Then $T^{4n}(0)=0$ and $T^{4n+1}(0)=c$ for all positive integers $n$, showing that $T^k(0),k=1,2,\cdots$ is neither convergent to 0 nor bounded away from zero.
\end{ex}
\begin{df}\rm Let $X$ be a Banach space and $A:X\to X$ a linear operator. We say $A$ has property (P) if
$X$ is a direct sum of two closed subspaces $V_1, V_2$ such that
(1) each $V_i$ is invariant under $A$, (2) $V_1$ is finite
dimensional, and (3) there exists $0\le r<1$ and a positive
integer $N$ such that $\|A^kx\|\le r^k\|x\|$ for all $x\in V_2$
and all $k\ge N$.
\end{df}
Note that by Gelfand's spectral radius theorem, condition (3)
above is equivalent to that  $A$, as an operator on $V_2$, has spectral
radius less than 1.\\
It is well-known that every compact operator, or more generally,
Riesz operator, has property (P), see e.g. [1].

\begin{theo} Let $(X,\|\cdot\|)$ be a Banach space and $A:X\to X$ a bounded
operator having property (P). Let $c$ be a constant vector in $X$,
and let $T(x)=Ax+c$ for $x\in X$. Then for any vector $x\in X$,
either $\{T^kx, k=0,1,\cdots\}$ is bounded or
$\lim_{k\to\infty}\|T^kx\|=\infty$. Moreover, if $c=0$, then either $\lim_{k\to\infty}T^kx=0$ or  $\lim_{k\to\infty} \|T^k x\|=\infty$ or $F\le\|T^kx\|\le G$ for sufficiently large $k$'s, where $0<F\le G<\infty$.
\end{theo}
{\bf Proof.}\\
Let $V_i,i=1,2, N, r$ be as in Definition 1. Let $P_i, i=1,2$ be the projections of $X$ onto $V_i, i=1,2$ respectively. Define a norm $|\cdot|$ on $X$ as \[|x|=\|P_1x\|+\|P_2x\|.\]
$|\cdot|$ is equivalent to $\|\cdot\|$ and $T$ commutes with $P_i,i=1,2$. \\
For any $v\in X$, write $v_i=P_iv, i=1,2$. Then for $x\in X$, and $k\ge N$ one has
\begin{eqnarray*} \|T^kx_2\| & = & \|A^kx_2+c_2+Ac_2+\cdots+A^{k-1}c_2\|\\
&\le& r^k\|x_2\|+ \|c_2+\cdots +A^{N-1}c_2\|+(r^N+\cdots+r^{k-1})\|c_2\|\\
&\le& r^k\|x_2\|+\|c_2+\cdots +A^{N-1}c_2\|+\frac{1}{1-r}\|c_2\|,
\end{eqnarray*}
showing that $T^kx_2,k=1,2,\cdots$ is bounded. By Theorem 1, $\|T^kx_1\|, k=1,2,\cdots$ is either bounded or approaching infinity. Hence $|T^kx|=\|T^kx_1\|+\|T^kx_2\|, k=1,2,\cdots$ is either bounded or approaching infinity. Since norms $|\cdot|$ and $\|\cdot\|$ are equivalent, the same is true for $\|T^kx\|, k=1,2,\cdots$.\\
If $c=0$, then $c_2=0$ and $T^kx_2\to 0$ as $k\to\infty$, and the last part of the theorem follows readily from Theorem 2.
\begin{co} Let $(X,\|\cdot\|)$ be a Banach space and $A:X\to X$ a Riesz
operator. Let $c$ be a constant vector in $X$,
and let $T(x)=Ax+c$ for $x\in X$. Then for any vector $x\in X$,
either $\{T^kx, k=0,1,\cdots\}$ is bounded or
$\lim_{k\to\infty}\|T^kx\|=\infty$. Moreover, if $c=0$, then either $\lim_{k\to\infty}T^kx=0$ or  $\lim_{k\to\infty} \|T^k x\|=\infty$ or $F\le\|T^kx\|\le G$ for sufficiently large $k$'s, where $0<F\le G<\infty$.
\end{co}
\begin{co} Let $(X,\|\cdot\|)$ be a Banach space and $A:X\to X$ a compact
operator. Let $c$ be a constant vector in $X$,
and let $T(x)=Ax+c$ for $x\in X$. Then for any vector $x\in X$,
either $\{T^kx, k=0,1,\cdots\}$ is bounded or
$\lim_{k\to\infty}\|T^kx\|=\infty$. Moreover, if $c=0$, then either $\lim_{k\to\infty}T^kx=0$ or  $\lim_{k\to\infty} \|T^k x\|=\infty$ or $F\le\|T^kx\|\le G$ for sufficiently large $k$'s, where $0<F\le G<\infty$.
\end{co}
Let us consider now the behavior of the sequence of averages of T:
\[\mbox{Ave}_kT(x)=\frac{1}{k}(x+Tx+\cdots+T^{k-1}x), k=1,2,\cdots \]We have
\[\mbox{Ave}_kT(x)=\frac{1}{k}(x+Ax+\cdots+A^{k-1}x)+\frac{1}{k}((k-1)c+(k-2)Ac+\cdots+A^{k-2}c)\]
As in the proof of Theorem 1, we may assume that $A=\lambda I+N$, and that $s,t$ are defined as in there. By replacing $c$ by $x$, $x$ by 0, and hence $s$ by 0 and $t$ by $s$ in the proof of Theorem 1, we obtain:\\
For $\lambda\neq 1$,\[\frac{1}{k}(x+Ax+\cdots+A^{k-1}x)=\frac{1}{k}B+\sum_{j=0}^{s-1}\lambda^{k-j}\frac{C(k,j)}{k}A_j\] where \[A_j=-\sum_{i=j}^{s-1}\frac{1}{(1-\lambda)^{i-j+1}}N^{i}x, j=0,\cdots,s-1\]
\[B=\sum_{j=0}^{s-1}\frac{1}{(1-\lambda)^{j+1}}N^{j}x.\]
For $\lambda=1$, \[\frac{1}{k}(x+Ax+\cdots+A^{k-1}x)=\frac{1}{k}\sum_{j=1}^sC(k,j)N^{j-1}x.\]
By expanding $A^j$ we have \[(k-1)c+(k-2)Ac+\cdots+A^{k-2}c=T_0c+T_1Nc+\cdots+T_{t-1}N^{t-1}c\] where
\[T_j=\sum_{i=j}^{k-2}(k-i-1)C(i,j)\lambda^{i-j}, j=0,\cdots,t-1.\]
Assume that $\lambda\neq 1$. By subtracting $\lambda T_0$ from $T_0$ and using the  geometric series formula,
one gets \[T_0=\frac{k}{1-\lambda}-\frac{1-\lambda^k}{(1-\lambda)^2}\]
Using the relation $(i-j)C(i,j)=(j+1)C(i,j+1)$, we see that
\[T_{j+1}=\frac{1}{j+1}\frac{d}{d\lambda}T_j\] We shall prove in the Appendix  that
\begin{equation}T_j=
\frac{k}{(1-\lambda)^{j+1}}-\frac{j+1}{(1-\lambda)^{j+2}}+\frac{C(k,j)\lambda^{k-j}}{(1-\lambda)^2}D(k,j,\lambda)\end{equation} for $j=0,\cdots,t-1$, where \[D(k,j,\lambda)=\frac{1}{(1-\lambda)^j}[B_0(k,j)\lambda^j+B_1(k,j)\lambda^{j-1}+\cdots+B_{j-1}(k,j)\lambda+1]\] and \[
B_i(k,j)=(-1)^{j-i}C(j,i)\frac{C(k-j,2)}{C(k-i,2)},\
i=0,\cdots,j\]
Note that \[B_i(k,j)\to (-1)^{j-i}C(j,i)\] so that $D(k,j,\lambda)$ approaches 1 as $k\to\infty$. Substituting the formulae we obtain thus far into $\mbox{Ave}_kT(x)$, with $w=\max\{s,t\}$, we get for $\lambda\neq 1$,
\begin{equation}\mbox{Ave}_kT(x)=E+\frac{1}{k}F+G(k,\lambda)\end{equation} where \begin{equation}G(k,\lambda)=\sum_{j=0}^{w-1}\frac{1}{k}C(k,j)\lambda^{k-j}\left(\epsilon(s-1-j)A_j
+\epsilon(t-1-j)\frac{D(k,j,\lambda)}{(1-\lambda)^2}N^jc\right)\end{equation}
\[E=\sum_{j=0}^{t-1}\frac{1}{(1-\lambda)^{j+1}}N^jc,\]
\[F=\sum_{j=0}^{s-1}\frac{1}{(1-\lambda)^{j+1}}N^jx-\sum_{j=0}^{t-1}\frac{j+1}{(1-\lambda)^{j+2}}N^jc,\]
and $\epsilon(r)=1$ for $r\ge 0$ and $\epsilon(r)=0$ for $r<0$. \\
If $\lambda=1$, then $T_j$ is equal to
\[S_{j,k}=C(k-2,j)+2C(k-3,j)+\cdots+(k-j-1)C(j,j)\]
We have \[S_{j,k}-S_{j,k-1}=C(k-2,j)+C(k-3,j)+\cdots+C(j,j)=C(k-1,j+1)\] where the last equality follows from repeatedly applying the identity \\ $C(j,i+1)+C(j,i)=C(j+1,i+1)$. From this it follows that
\begin{eqnarray*}S_{j,k} &=& S_{j,k-1}+C(k-1,j+1)\\
&=&S_{j,k-2}+C(k-2,j+1)+C(k-1,j+1)\\
&=& \cdots\\
&=& S_{j,j+2}+C(j+2,j+1)+\cdots+C(k-1,j+1)\\
&=& C(j+1,j+1)+C(j+2,j+1)+\cdots+C(k-1,j+1)\\
&=& C(k,j+2)
\end{eqnarray*}
Therefore  for $\lambda=1$,
\begin{eqnarray}\mbox{Ave}_kT(x)&=&\frac{1}{k}[C(k,1)x+C(k,2)Nx+\cdots+C(k,s)N^{s-1}x]\\& &+\frac{1}{k}[C(k,2)c+C(k,3)Nc+\cdots+C(k,t+1)N^{t-1}c]\end{eqnarray}

Case 1: $|\lambda|>1$. Suppose that $t>s$. Then from (7) and (8),
$\frac{1}{k}\lambda^kC(k,t-1)$ dominates all other coefficients.
Since $D(k,t-1,\lambda)\to 1$ as $k\to\infty$ and $N^{t-1}c\neq
0$, we see that $\|\mbox{Ave}_kT(x)\|\to \infty$ as $k\to\infty$.
The same is true if $s>t$ since $A_{s-1}\neq 0$. So assume that
$s=t$. In the trivial case $s=t=0$, we have $\mbox{Ave}_kT(x)=0$
for all $k$. So assume that $s=t\ge 1$. Direct checking (see Appendix, item 1) shows that $C(k,j)C(j,i)\frac{C(k-j,2)}{C(k-i,2)}=C(k,i)C(k-i-2,j-i)$ so that $C(k,j)D(k,j,\lambda)$ is a polynomial in $k$. Then from (8), we have
\begin{eqnarray*}G(k,\lambda)&=&\sum_{j=0}^{t-1}\frac{1}{k}C(k,j)\lambda^{k-j}\left(A_j
+\frac{D(k,j,\lambda)}{(1-\lambda)^2}N^jc\right)\\
&=&\frac{\lambda^k}{k}\sum_{j=0}^{t-1}p_1(k,j,\lambda)
A_j+p_2(k,j,\lambda)N^jc
 \end{eqnarray*} where for fixed $\lambda$,
$p_1,p_2$ are polynomials in $k$.
Since $\lim_{k\to\infty}\frac{\lambda^k}{k}p(k)=\infty$ for any polynomial $p(k)$ we see that if
\[H(k,\lambda)=\sum_{j=0}^{t-1}p_1(k,j,\lambda) A_j+p_2(k,j,\lambda)N^jc,\] as a polynomial in $k$, is
identically zero, then clearly from (7) we have $\mbox{Ave}_kT(x)\to E$, otherwise $\|\mbox{Ave}_kT(x)\|\to
\infty$.  \\

Case 2: $|\lambda|=1, \lambda\neq 1$. Suppose that $t>s$.
Since $D(k,t-1,\lambda)\to 1$ as $k\to\infty$ and $N^{t-1}c\neq
0$, we see that $\|\mbox{Ave}_kT(x)\|\to \infty$ as $k\to\infty$ if $t\ge 3$, and it approaches to
0 or is bounded if $t\le 2$.
The same is true if $s>t$ since $A_{s-1}\neq 0$. So assume that
$s=t$. In the trivial case $s=t=0$, we have $\mbox{Ave}_kT(x)=0$
for all $k$. So assume that $s=t\ge 1$. As in Case 1, we consider the polynomial (in $k$) \[H(k,\lambda)=\sum_{j=0}^{t-1}C(k,j)\lambda^{-j}\left(A_j
+\frac{D(k,j,\lambda)}{(1-\lambda)^2}N^jc\right)\] so that
\[G(k,\lambda)=\frac{\lambda^k}{k}H(k,\lambda)\]
If degree of $H$ is one or less, then we see from above that $G(k,\lambda)$ is bounded and hence $\mbox{Ave}_kT(x)$ is bounded by (7). If degree of $H$ is two or more, then $\|G(k,\lambda)\|\to \infty$ as $k\to\infty$ and hence so is $\|\mbox{Ave}_kT(x)\|$. \\

Case 3: $\lambda=1$. We refer to (9) and (10). If $s=t=0$, then $c=x=0$ and  $\mbox{Ave}_kT(x)=0$. If $t=0$ and $s=1$, then $\mbox{Ave}_kT(x)=x$. If $t=s-1, s\ge 2$ and $N^ic+N^{i+1}x=0$ for all $i=0,\cdots,t-1$, then $\mbox{Ave}_kT(x)=x$. In all other cases $\lim_k\|\mbox{Ave}_kT(x)\|=\infty$.\\

The previous discussions yield the proof of the following:
\begin{theo} Let $T: \mathbb{C}^n \to \mathbb{C}^n$ be an affine map defined by $Tx=Ax+c$, where $A$ is an $n\times n$ complex matrix and $c$ a constant vector in $\mathbb{C}^n$. Let $\|\cdot\|$ be a norm on $\mathbb{C}^n$. For any vector $x\in\mathbb{C}^n$, define \[\mbox{Ave}_kT(x)=\frac{1}{k}(x+Tx+\cdots + T^{k-1}x)\]
Then the sequence \[\mbox{Ave}_kT(x),k=0,1,2,\cdots\] is either bounded or $\lim_k \|\mbox{Ave}_kT(x)\|=\infty$.
\end{theo}
If $T$ is linear, i.e. if $c=0$, we can say more:
\begin{theo} Let $A: \mathbb{C}^n \to \mathbb{C}^n$ be linear map. Let $\|\cdot\|$ be a norm on $\mathbb{C}^n$. For any vector $x\in\mathbb{C}^n$, define \[\mbox{Ave}_kA(x)=\frac{1}{k}(x+Ax+\cdots + A^{k-1}x)\]
Then the sequence \[\mbox{Ave}_kA(x),k=0,1,2,\cdots\] is either (i) convergent to 0, or (ii) $\lim_k \|\mbox{Ave}_kA(x)\|=\infty$, or (iii) $F\le\|\mbox{Ave}_kAx\|\le G$ for sufficiently large $k$'s, where $F,G$ are positive numbers with $ F\le G$.
\end{theo}{\bf Proof.}\\
If $c=0$ in Theorem 4, then $E=0$, and by examining its proof we see that   in all cases where $T^kx$ is bounded and not convergent to 0, it is bounded away from 0. (In Case 2 of the proof, if $H(k,\lambda)$ is of degree one, then $G(k,\lambda)$ is bounded away from 0, and if $H(k,\lambda)$ is of degree 0, i.e. a constant vector, then $G(k,\lambda)\to 0$.)  Q.E.D.\\ \\

The following example shows that Theorem 5 is false if the map is not linear.\\
\begin{ex}\rm Let \[A=\left(\begin{array}{cc}i & 1\\0&i\end{array}\right)\] and let $c=\left(\begin{array}{c}1\\ \\0\end{array}\right)$ Define $T:\mathbb{C}^2 \to \mathbb{C}^2$ by $Tx=Ax+c$. Let $v=\left(\begin{array}{c}0\\ \\1\end{array}\right)$. Then
\[\mbox{Ave}_kT(v)=\frac{1}{k}\frac{1-i^k}{1-i}v+\frac{1-i^{k-1}}{1-i}\]
which does not converge to 0, and has a subsequence converging to 0.
\end{ex}
Theorems 4 and 5 are also valid in Banach spaces. For proof, one only has to use these theorems and note that the average of the bounded (resp. null convergent) sequence $T^kx_2,k=0,1,\cdots$ in the proof of Theorem 3 is also bounded (resp. null convergent). Thus we have
\begin{theo} Let $X$ be a Banach space. Let $T: X \to X$ be an affine map defined by $Tx=Ax+c$, where $A$ is an operator with property (P) and $c$ a constant vector in $X$. Let $\|\cdot\|$ be the norm on $X$. For any vector $x\in X$, define \[\mbox{Ave}_kT(x)=\frac{1}{k}(x+Tx+\cdots + T^{k-1}x)\]
Then the sequence \[\mbox{Ave}_kT(x),k=0,1,2,\cdots\] is either bounded or $\lim_k \|\mbox{Ave}_kT(x)\|=\infty$. Moreover, if $c=0$, then  the sequence \[\mbox{Ave}_kA(x),k=0,1,2,\cdots\] is either (i) convergent to 0, or (ii) $\lim_k \|\mbox{Ave}_kA(x)\|=\infty$, or (iii) $F\le\|\mbox{Ave}_kAx\|\le G$ for sufficiently large $k$'s, where $0<F\le G<\infty.$
\end{theo}
Recall that compact operators, or more generally, Riesz operators have property (P). So the above theorem is valid for these operators. \\

Our last objective is to prove the following theorem. The result concerns $\mbox{Ave}_kT(x)$ in the case $s=t$. It shows in particular that $H$ is identically 0 if and only if
$x=E$, the unique fixed point of $T$, i.e. $Tx=x$, so that  in case (i) in the proof of Theorem 4  we actually have $\mbox{Ave}_kT(x)= E$ for all $k$, not just $\mbox{Ave}_kT(x)\to E$ as $k\to \infty$.\\
\begin{theo} Let $\lambda$ be a fixed complex number, $\lambda\neq 0,1$. Let \[H(k)=\sum_{j=0}^{t-1}C(k,j)\lambda^{-j}\left(A_j
+\frac{D(k,j,\lambda)}{(1-\lambda)^2}N^jc\right)\]
 be defined as previously.  Let $1\le i\le t-1$. The degree of $H(k)$ is at most $i-1$ if and only if \begin{equation}N^ix=\frac{1}{(1-\lambda)}N^ic+ \frac{1}{(1-\lambda)^2}N^{i+1}c+\cdots+\frac{1}{(1-\lambda)^{t-i}}N^{t-1}c\end{equation}
$H(k)$ is identically 0 if and only if
\begin{equation}x=\frac{1}{(1-\lambda)}c+ \frac{1}{(1-\lambda)^2}Nc+\cdots+\frac{1}{(1-\lambda)^{t}}N^{t-1}c\end{equation}
\end{theo}
The proof will follow from the following discussions.\\ \\
Consider $D(k,j,\lambda)-1$. We have
\[D(k,j,\lambda)-1=\frac{1}{(1-\lambda)^j}\sum_{i=0}^j(-1)^iC(j,i)\left(\frac{C(k-j,2)}{C(k-j+i,2)}-1\right)\lambda^i\]
and for $k>j$,
\begin{eqnarray*}& &\frac{C(k-j,2)}{C(k-j+i,2)}-1\\&=&i\frac{-2k+2j-i+1}{(k-j+i)(k-j+i-1)}\\
&=&i\left(\frac{-2}{k}+\frac{2j-i+1}{k^2}\right)\left(1-\frac{j-i}{k}\right)^{-1}\left(1-\frac{j-i+1}{k}\right)^{-1}\\
&=& \frac{i}{k}\left(-2+\frac{2j-i+1}{k}\right)\left(1+\frac{j-i}{k}+\frac{(j-i)^2}{k^2}+\cdots\right)\left(1+\frac{j-i+1}{k}+\frac{(j-i+1)^2}{k^2}+\cdots\right)\\
&=& i\sum_{m=1}^\infty \frac {-2(j-i)^{m-1}-(i-1)((j-i)^{m-1}-(j-i+1)^{m-1})}{k^m}
\end{eqnarray*}
The coefficient of $\frac{1}{k^m}$ in the sum above is obtained from the following calculations:\\
\begin{eqnarray*}& & -2\sum_{p=0}^{m-1} (j-i)^p(j-i+1)^{m-1-p}+(2j-i+1)\sum_{p=0}^{m-2}(j-i)^p(j-i+1)^{m-2-p}\\
&=& -2(j-i)^{m-1}-2(j-i+1)\sum_{p=0}^{m-2} (j-i)^p(j-i+1)^{m-2-p}+(2j-i+1)\sum_{p=0}^{m-2}(j-i)^p(j-i+1)^{m-2-p}\\
&=& -2(j-i)^{m-1}+(i-1)\sum_{p=0}^{m-2}(j-i)^p(j-i+1)^{m-2-p}\\
&=& -2(j-i)^{m-1}-(i-1)((j-i)^{m-1}-(j-i+1)^{m-1})
\end{eqnarray*}
Our first objective is to write $(j-i)^{m-1}$ as a combination of $x_0=1, x_1=i-1, x_2=(i-2)(i-1),\cdots, x_{m-1}=(i-m+1)\cdots(i-1)$. The proof for the following item 1 can be found in [3].
\begin{enumerate}
\item Let $f(x)=a_lx^l+\cdots+a_1x+a_0$ be a polynomial of degree $l$. Let $x_0=1, x_1=x-1, x_2=(x-2)(x-1),\cdots, x_l=(x-l)\cdots(x-1)$.Then
\[f(x)=d_lx_l+\cdots +d_0x_0\] where
\begin{eqnarray*}d_0 &=& f(1)\\
d_1 &=& f(2)-f(1)\\
d_2 &=& \frac{1}{2!}\left[f(3)-2f(2)+f(1)\right]\\
d_3 &=& \frac{1}{3!}\left[f(4)-3f(3)+3f(2)-f(1)\right]\\
& &\cdots\cdots\cdots\cdots\cdots\cdots\\
d_l &=&\frac{1}{l!}\left[f(l+1)-C(l,1)f(l)+\cdots+(-1)^lf(1)\right]\\
&=& a_l
\end{eqnarray*}
Moreover, for any integer $n$ with $n>l$, \[f(n+1)-C(n,1)f(n)+C(n,2)f(n-1)+\cdots+(-1)^nf(1)=0\]

\item Applying above to the polynomial $f(i)=(j-i)^{m-1}$, we get \[(j-i)^{m-1}=c_{m-1}x_{m-1}+\cdots +c_0x_0\] where
\begin{eqnarray*}c_0 &=& (j-1)^{m-1}\\
c_1 &=& (j-2)^{m-1}-(j-1)^{m-1}\\
c_2 &=& \frac{1}{2!}\left[(j-3)^{m-1}-2(j-2)^{m-1}+(j-1)^{m-1}\right]\\
c_3 &=& \frac{1}{3!}\left[(j-4)^{m-1}-3(j-3)^{m-1}+3(j-2)^{m-1}-(j-1)^{m-1}\right]\\
& &\cdots\cdots\cdots\cdots\cdots\cdots\\
c_{m-1} &=&\frac{1}{(m-1)!}\left[(j-m)^{m-1}-C(m-1,1)(j-m+1)^{m-1}+\cdots+(-1)^{m-1}(j-1)^{m-1}\right]\\
&=& (-1)^{m-1}
\end{eqnarray*}
Moreover, for any $n>m-1$, \[(j-n-1)^{m-1}-C(n,1)(j-n)^{m-1}+C(n,2)(j-n+1)^{m-1}+\cdots+(-1)^n(j-1)^{m-1}=0\]

\item Let $x_i,c_i, i=0,\cdots,m-1$ be defined above. Then
\[p(i)=(i-1)\left[(j-i)^{m-1}-(j-i+1)^{m-1}\right]=c_1x_1+2c_2x_2+\cdots+(m-1)c_{m-1}x_{m-1}\]
so that
\[-2(j-i)^{m-1}-p(i)=-2c_0x_0-3c_1x_1-\cdots-(m+1)c_{m-1}x_{m-1}\]
{\bf Proof.}\\
Note that the degree of $p(i)$ is $m-1$ and the coefficient of $i^{m-1}$ is $(-1)^{m-1}(m-1)$. Applying item 1 to $p(i)$ and using the identity $C(k-1,i)+C(k-1,i+1)=C(k,i+1)$ we get
\begin{eqnarray*}
d_0 &=& p(1)=0\\
d_1 &=& p(2)-p(1)=(j-2)^{m-1}-(j-1)^{m-1}=c_1\\
d_2 &=& \frac{1}{2!}\left[p(3)-2p(2)+p(1)\right]=(j-3)^{m-1}-(j-2)^{m-1}-[(j-2)^{m-1}-(j-1)^{m-1}]=2c_2\\
& &\cdots\cdots\cdots\cdots\cdots\cdots\\
d_k &=& \frac{1}{k!}\left[p(k+1)-C(k,1)p(k)+\cdots+(-1)^iC(k,i)p(k-i+1)+\cdots+(-1)^{k-1}kp(2)+\right.\\
& &\left.(-1)^kp(1)\right]\\
&=&\frac{1}{k!}k\left[(j-k-1)^{m-1}-C(k,1)(j-k)^{m-1}+\cdots+(-1)^k(j-1)^{m-1}\right]=kc_k\\
& &\cdots\cdots\cdots\cdots\cdots\cdots\\
d_{m-1} &=& (-1)^{m-1}(m-1)=(m-1)c_{m-1}
\end{eqnarray*}
\item Write $i_0=i, i_1=i(i-1),\cdots, i_p=i(i-1)\cdots(i-p)$. Note that $i_p$ is $i$ times $x_p$ in item 2. This change is necessary because there is an $i$ that was factored out of the summation sign in our expansion of  $\frac{C(k-j,2)}{C(k-j+i,2)}-1$.\\
By item 3, the coefficient of $\frac{1}{k^m}$ in the expansion of $\frac{C(k-j,2)}{C(k-j+i,2)}-1$ is
$-\sum_{p=0}^{m-1}(p+2)c_pi_p$. Therefore the coefficient $L(m,j,\lambda)$ of $\frac{1}{k^m}$ in the expansion of $D(k,j,\lambda)-1$ is given by \begin{eqnarray*}
    & &-\frac{1}{(1-\lambda)^j}\sum_{i=0}^j(-1)^iC(j,i)\left(\sum_{p=0}^{m-1}(p+2)c_pi_p\right)\lambda^i\\
    &=&-\frac{1}{(1-\lambda)^j}\sum_{p=0}^{m-1}(p+2)c_p\lambda^{p+1}\sum_{i=0}^j(-1)^iC(j,i)i_p\lambda^{i-p-1}\\
    &=& -\frac{1}{(1-\lambda)^j}\sum_{p=0}^{m-1}(p+2)c_p\lambda^{p+1}\frac{d^{p+1}}{d\lambda^{p+1}}(1-\lambda)^j\\
    &=&\frac{1}{(1-\lambda)^j}\sum_{p=0}^{m-1}(-1)^p(p+2)c_p\lambda^{p+1}j(j-1)\cdots(j-p)(1-\lambda)^{j-p-1}\\
    &=&\sum_{p=0}^{m-1}(-1)^p(p+2)c_pj(j-1)\cdots(j-p)\left(\frac{\lambda}{1-\lambda}\right)^{p+1}\\
    &=&\sum_{i=1}^m (-1)^{i+1}(i+1)!c_{i-1}C(j,i)\left(\frac{\lambda}{1-\lambda}\right)^i
    \end{eqnarray*}
To emphasize that $c_{i-1}$ depends on $j,m$, we write $c_i$ as $P(i,j,m)$, so
\[P(i,j,m)=  \frac{1}{i!}[(j-i-1)^{m-1}-C(m-1,1)(j-i)^{m-1}+\cdots+(-1)^{m-1}(j-1)^{m-1}] \]
and \[L(m,j,\lambda)=\sum_{i=1}^m (-1)^{i+1}(i+1)!P(i-1,j,m)C(j,i)\left(\frac{\lambda}{1-\lambda}\right)^i\]
We have $P(0,j,m)=c_0=(j-1)^{m-1}, P(m-1,j,m)=c_{m-1}=(-1)^{m-1}, P(i,j,m)=c_i=0$ for $i\ge m$. We define $P(i,j,m)=0$ for negative $i$. We have the following recursive relations, which will not be used in the sequel,
\[P(i,j,m)=(j-i-1)P(i,j,m-1)-P(i-1,j,m-1), m\ge 2.\]
Also  we define
\[M(m,j,x)=\sum_{i=1}^m (-1)^{i+1}(i+1)!P(i-1,j,m)C(j,i)x^i\]
so that \[L(m,j,\lambda)=M(m,j,\frac{\lambda}{1-\lambda})\]
\item
\begin{eqnarray*}& &M(m,j,x)-H_1(j_1,\cdots,j_{m-1})M(m-1,j,x)+\cdots\\&+&(-1)^iH_i(j_1,\cdots,j_{m-1})M(m-i,j,x)+\cdots +(-1)^{m-1}H_{m-1}(j_1,\cdots,j_{m-1})M(1,j,x)\\&=&j_0j_1\cdots j_{m-1}(m+1)x^m\end{eqnarray*}
where $j_i=j-i,i=0,\cdots,m-1$ and $H_i(j_1,\cdots,j_{m-1})$ is the unsigned coefficient of $y^i$ in the expansion of
the polynomial
\[(y-j_1)(y-j_2)\cdots(y-j_{m-1})\] i.e. $H_1=j_1+\cdots+j_{m-1}, H_2=\sum_{i<k}j_ij_k, \cdots, H_{m-1}=j_1\cdots j_{m-1}$.\\
In particular, we have
\begin{eqnarray*}
M(1,j,x)&=&2jx\\
M(2,j,x)&=&j(j-1)3x^2+(j-1)M(1,j,x)=j(j-1)(3x^2+2x)\\
M(3,j,x)&=&j(j-1)(j-2)4x^3+(j_1+j_2)M(2,j,x)-j_1j_2M(1,j,x)\\&=&j(j-1)((j-2)4x^3+(2j-3)3x^2+2(j-1)x)
\end{eqnarray*}
{\bf Proof.}\\
For fixed $m,j$, \begin{eqnarray*}p(x)&=&M(m,j,x)-H_1(j_1,\cdots,j_{m-1})M(m-1,j,x)+\cdots\\&+&(-1)^iH_i(j_1,\cdots,j_{m-1})M(m-i,j,x)+\cdots +(-1)^{m-1}H_{m-1}(j_1,\cdots,j_{m-1})M(1,j,x)\end{eqnarray*} is a polynomial in $x$ with $p(0)=0$. Fix $1\le k\le m-1$. We shall show that the coefficient of $x^k$ in $p(x)$ is 0. Now the coefficient is \\ $(-1)^{k+1}(k+1)!C(j,k)a_k$, where
\begin{eqnarray*}a_k&=&P(k-1,j,m)-H_1P(k-1,j,m-1)+\cdots+(-1)^iH_iP(k-1,j,m-i)+\cdots\\&+&(-1)^{m-1}H_{m-1}P(k-1,j,1)\end{eqnarray*}
Since \[y^{m-1}-H_1y^{m-2}+\cdots+(-1)^iH_iy^{m-i-1}+\cdots+(-1)^{m-1}H_{m-1}=(y-j_1)\cdots(y-j_{m-1})\] we have for each $1\le z\le k$,\[j_z^{m-1}-H_1j_z^{m-2}+\cdots+(-1)^iH_ij_z^{m-i-1}+\cdots+(-1)^{m-1}H_{m-1}=0\]
From the  definition of $P(i,j,m)$, $j_z^{m-1}, j_z^{m-2},\cdots$ are the corresponding terms (with the same coefficient) in $P(k-1,j,m),P(k-1,j,m-1),\cdots$. It follows that $a_k=0$. \\
For $k=m$, we have $a_k=P(m-1,j,m)=(-1)^{m-1}$ since $P(i,j,m)=0$, for $i\ge m$. Thus the coefficient of $x^m$ is $(-1)^{m+1}(-1)^{m-1}(m+1)!C(j,m)=(m+1)j_0j_1\cdots j_{m-1}$. This completes the proof.

\item Let $m_1,m_2,\cdots,$ be variables. Let $k,j$, be positive integers. Define $R(k,j)$ to be the sum of \[m_1^am_2^b\cdots m_k^c\] where $a,b,\cdots,c$ are nonnegative integers such that $0\le a,b,\cdots,c\le j$ and $a+b+\cdots +c=j$. Define $R(k,0)=1$. Define $S(k,j)$ to be the sum of \[m_1^am_2^b\cdots m_k^c\] where $a,b,\cdots,c$ are nonnegative integers such that $0\le a,b,\cdots,c\le 1$ and $a+b+\cdots +c=j$. Define $S(k,0)=1$ and $S(0,0)=1$. Then we have the following identity:
    \[\sum_{i=0}^q (-1)^iR(p+i,q-i)S(p+i-1,i)=0\]
    for any positive integers $p,q$.\\
{\bf Proof.}\\
Consider a typical term $t=m_1^am_2^b\cdots m_{p+i}^k$ resulting from the summand
$(-1)^iR(p+i,q-i)S(p+i-1,i)$, sign disregarded. Let us write $c(1)=a, c(2)=b,\cdots, c(p+i)=k$. We may assume that $c(p+i)\ge 1$, otherwise the term belongs to a previous summand; this assumption also implies that $t$ does not belong to any previous (smaller $i$) summand. Let
\[ A=\{j: 1\le j\le p+i-1, c(j)\neq 0\}\]
Then the cardinality $|A|$ of $A$ must be greater than or equal to $i$ because of $S(p+i-1,i)$. Also the term $t$ appears in the expansion of the summand for exactly $C(|A|,i)$ times. Next let us consider how many times $t$ appears in the expansion of the next summand $(-1)^{i+1} R(p+i+1, q-i-1)S(p+i,i+1)$. $t$ can result from multiplying a term in $R(p+i+1,q-i-1)$ with a term $s$ in $S(p+i,i+1)$. Denote the exponent of $m_{p+i}$ in $s$ by $b_s(p+i)$, which is either 1 or 0. \\
Denote by $S^*$ and $T^*$ the set of terms in $S(p+i,i+1)$ and $R(p+i+1,q-i-1)$ respectively. Consider the following two sets:
\[U_k=\{s: s\in S^*, b_s(p+i)=k, rs=t \mbox{ for some } r\in R^*\}, k=0,1.\]
Note that  each $s\in U_k$ corresponds to exactly one $r$ such that $rs=t$. Clearly $U_1$ has exactly $C(|A|,i)$ elements, and these elements yields the same number of $t$'s which cancel out with the previous ones because of the sign change.\\
Each $s\in U_0$ consists of $i+1$ factors from $m_1,\cdots, m_{p+i-1}$, so $U_0$ has $C(|A|,i+1)$ elements which yield the same number of $t$'s.\\
If $C(|A|,i+1)=0$, i.e. $|A|<i+1$, $S_0$ is empty and we are done since no further $t$'s will result. If not, we consider the next summand $(-1)^{i}R(p+i+2, q-i-2)S(p+i+1,i+2)$. \\
Denote by $S_1^*$ and $T_1^*$ the set of terms in $S(p+i+1,i+2)$ and $R(p+i+2,q-i-2)$ respectively. Consider the following two sets:
\[V_k=\{s: s\in S_1^*, b_s(p+i)=k, rs=t \mbox{ for some } r\in R_1^*\}, k=0,1.\]
Note that it is $b_s(p+i)$, not $b_s(p+i+1)$, in the above definition of $V_k$. Also note that for each $s\in V_k$, $b_s(p+i+1)$ must be 0 or else $rs=t$ is impossible. Then it is clear that $|V_1|=C(|A|,i+1)$ and $|V_0|=C(|A|,i+2),$ as before. The $C(|A|,i+1)$ $t$'s resulting from $V_1$ cancel out with the previous the same number of $t$'s. The same statement about $U_0$ above applies to $V_0$ and the process continues.
This process will continue for at most $p-1$ times since $|A|\le i+p-1$. This proves that the $t$'s occurring in $\sum_{i=0}^q (-1)^iR(p+i,q-i)S(p+i-1,i)$ are all canceled out. Q.E.D.
\item We now finish the proof of Theorem 7. Let $i=t-1$. $H$ has degree of at most $i-1$ if and only if
the coefficient of $k^{t-1}$ is 0. Since $\lim_{k\to\infty}D(k,j,\lambda)=1$, this amounts to \[A_{t-1}+\frac{1}{(1-\lambda)^2}N^{t-1}c=0\] which is equivalent to \[N^{t-1}x=\frac{1}{1-\lambda}N^{t-1}c\]
So the assertion in Theorem 7 is true for $i=t-1$. Suppose the assertion is true for some $i\ge 1$; this implies that coefficients of $k^j$ in $H$ are 0 for $j\ge i$.  We shall prove that it is also true for $i-1$. This will complete the proof by induction.\\
By induction hypothesis, we have
\[N^ix=\sum_{j=i}^{t-1}\frac{1}{(1-\lambda)^{j-i+1}}N^jc\]
By applying $N$ to both sides $t-1-i$ times, recalling that $N^tc=0$, we get
\begin{equation}N^l x= \sum_{j=l}^{t-1} \frac{1}{(1-\lambda)^{j-l+1}}N^jc\end{equation}
for all $l=i,\cdots,t-1$. \\
Substituting (13) into the formula for $A_l$, we find
\begin{equation}A_l=-\sum_{j=l}^{t-1}\frac{j-l+1}{(1-\lambda)^{j-l+2}}N^jc\end{equation}
for all $l=i,\cdots,t-1$. \\
Using $A_{j-1}=\frac{1}{1-\lambda}(A_j-N^{j-1}x)$, we get
\begin{equation}A_{i-1}=-\sum_{j=i}^{t-1}\frac{j-i+1}{(1-\lambda)^{j-i+3}}N^jc-\frac{1}{1-\lambda}N^{i-1}x\end{equation}
If we write $P_j$ for the term \[C(k,j)\lambda^{-j}\left(A_j
+\frac{D(k,j,\lambda)}{(1-\lambda)^2}N^jc\right)\] in $H$ and substitute (14) and (15), we find that, for $i\le l\le t-1$
\begin{equation}P_l=C(k,l)\lambda^{-l}\left[-\sum_{p=l+1}^{t-1}\frac{p-l+1}{(1-\lambda)^{p-l+2}}N^pc+\frac{D(k,l,\lambda)-1}{(1-\lambda)^2}N^lc\right]\end{equation}
where for $l=t-1$ the sum $\sum_{j=l+1}^{t-1}$ is an empty sum, and hence its value is 0; and  \begin{equation}P_{i-1}=C(k,i-1)\lambda^{-i+1}\left[-\sum_{j=i}^{t-1}\frac{j-i+1}{(1-\lambda)^{j-i+3}}N^jc-
\frac{1}{1-\lambda}N^{i-1}x+\frac{D(k,i-1,\lambda)}{(1-\lambda)^2}N^{i-1}c\right]\end{equation}
Since $C(k,j)$ and $C(k,j)D(k,j,\lambda)$ are polynomials in $k$ of degree $j$, the terms $P_j, j=0,\cdots,i-2$ in $H$ contain only $k$ powers of power less than $i-1$. And since the coefficients of $k$ powers of power greater than $i-1$ are 0 by induction hypothesis, we see that the coefficient of $k^{i-1}$ in $H$ is equal to the limit
\begin{equation}\lim_{k\to\infty} \frac{1}{k^{i-1}}(P_{i-1}+P_i+\cdots +P_{t-1})\end{equation}
Fix an $m$, $i\le m\le t-1$. Write $A=A(k,m)$ for the number
\[\frac{C(k,m)\lambda^{-m}}{k^m(1-\lambda)^2}\]
The coefficient of the vector $N^mc$ in $\frac{1}{k^{i-1}}P_j$ for each $j, i\le j\le m-1$, is, by (16), \[-C(k,j)\lambda^{-j}\frac{m-j+1}{(1-\lambda)^{m-j+2}}\] which can be rewritten as
\[A\frac{k^{m-j}}{(k-j)\cdots (k-m+1)}k^{j-i+1}(-(m-j+1)m(m-1)\cdots (j+1)(\frac{\lambda}{1-\lambda})^{m-j})\]
 and which by item 5 is equal to
\begin{eqnarray*}& &A\frac{k^{m-j}}{(k-j)\cdots (k-m+1)}k^{j-i+1}\left(-L(m-j,m,\lambda)+S(m-j-1,1)L(m-j-1,m,\lambda)+\cdots + \right.\\
& &\left.(-1)^{m-j}S(m-j-1,m-j-1)L(1,m,\lambda)\right)\\
&=&A\frac{k^{m-j}}{(k-j)\cdots (k-m+1)}k^{j-i+1}\sum_{l=1}^{m-j}(-1)^{m-j-l+1}S(m-j-1,m-j-l)L(l,m,\lambda)
\end{eqnarray*}
where $S$ is defined as in item 6, with $m_1=m-1,\cdots$. We also have
\begin{eqnarray*}\frac{k^{m-j}}{(k-j)\cdots (k-m+1)}k^{j-i+1}&=&k^{j-i+1}(1+R(m-j,1)k^{-1}+R(m-j,2)k^{-2}+\cdots\\& & +R(m-j,j-i+1)k^{-j+i-1})+O(1/k)\\
&=& \sum_{z=0}^{j-i+1}R(m-j,j-i+1-z)k^z +O(1/k)
\end{eqnarray*}
where $R$ is defined as in item 6, with $m_1=m-1,\cdots$.\\
The coefficient of the vector $N^mc$ in $\frac{1}{k^{i-1}}P_m$ is
\begin{eqnarray*}& &A k^{m-i+1}(D(k,m,\lambda)-1)\\
&=& Ak^{m-i}(L(1,m,\lambda)+L(2,m,\lambda)k^{-1}+\cdots +L(m-i+1,m,\lambda)k^{-m+i})+O(1/k)\\
&=& A \sum_{z=0}^{m-i}L(m-i+1-z,m,\lambda)k^z +O(1/k)
\end{eqnarray*}
Note that $N^mc$ does not appear in $P_l$ for $l>m$; so the coefficient of $N^mc$ in the sum \[\frac{1}{k^{i-1}}(P_i+\cdots +P_{t-1})\] is the same as that in the sum \[\frac{1}{k^{i-1}}(P_i+\cdots +P_{m}).\]
Therefore the coefficient of $N^mc$ in the sum \[\frac{1}{k^{i-1}}(P_i+\cdots +P_{t-1})\] which we shall call $c_m$ is
$A$ times \begin{eqnarray*}& &Y=\sum_{j=i}^{m-1}\sum_{l=1}^{m-j}(-1)^{m-j-l+1}S(m-j-1,m-j-l)L(l,m,\lambda)\sum_{z=0}^{j-i+1}R(m-j,j-i+1-z)k^z+\\& &\sum_{z=0}^{m-i}L(m-i+1-z,m,\lambda)k^z +O(1/k)\end{eqnarray*}
(Note that $A=O(1)$ so $AO(1/k)=O(1/k)$.)\\
Now we shall show that the polynomial part of $Y$ is a constant, i.e. the coefficients of $k^z, z=1,\cdots,m-i$ are zero. \\
Fix $z, 1\le z\le m-i$. Since $k^z$ occurs in the sum $\sum_{z=0}^{j-i+1}R(m-j,j-i+1-z)k^z$ only when $j-i+1\ge z$,i.e. $j\ge z+i-1$, we see that the coefficient of $k^z$ in $Y$ is $L(m-i+1-z,m,\lambda)$ plus \begin{eqnarray*}& & \sum_{j=z+i-1}^{m-1}\sum_{l=1}^{m-j}(-1)^{m-j-l+1}S(m-j-1,m-j-l)L(l,m,\lambda)R(m-j,j-i+1-z)\\
&=& \sum_{l=1}^{m-i-z+1}L(l,m,\lambda)\sum_{j=z+i-1}^{m-l}(-1)^{m-j-l+1}S(m-j-1,m-j-l)R(m-j,j-i+1-z)\\
&=&\sum_{l=1}^{m-i-z+1}L(l,m,\lambda)\sum_{n=0}^{m-i-z+1-l}(-1)^{n+1}R(l+n,m-i-z+1-l-n)S(l+n-1,n)\\
&=&\sum_{l=1}^{m-i-z}L(l,m,\lambda)\cdot 0+L(m-i-z+1,m,\lambda)\cdot (-1)\\
&=& -L(m-i-z+1,m,\lambda)
\end{eqnarray*}
where we have used item 5 for the 0. Hence the coefficient of $k^z$ is 0.\\
The constant term in $Y$ is given by the coefficient of $k^0$, which is \\$L(m-i+1,m,\lambda)$ plus
\begin{eqnarray*}
& &\sum_{j=i}^{m-1}\sum_{l=1}^{m-j}(-1)^{m-j-l+1}S(m-j-1,m-j-l)L(l,m,\lambda)R(m-j,j-i+1)\\
&=&\sum_{l=1}^{m-i}L(l,m,\lambda)\sum_{j=i}^{m-l}(-1)^{m-j-l+1}R(m-j,j-i+1)S(m-j-1,m-j-l)\\
&=&\sum_{l=1}^{m-i}L(l,m,\lambda)\sum_{n=0}^{m-i-l}(-1)^{n+1}R(l+n,m-l-i+1-n)S(l+n-1,n)\\
&=&\sum_{l=1}^{m-i}L(l,m,\lambda)(-1)^{m-i-l-1}R(l+n,0)S(m-i,m-i-l+1)\\
&=& \sum_{l=1}^{m-i}(-1)^{m-i-l-1}S(m-i,m-i-l+1)L(l,m,\lambda)
\end{eqnarray*}
 Adding the term $L(m-i+1,m,\lambda)$ back in we conclude that the said coefficient $c_m$ is
\[A[L(m-i+1,m,\lambda)-S(m-i,1)L(m-i,m,\lambda)+\cdots +(-1)^{m-i}S(m-i,m-i)L(1,m,\lambda)]+O(1/k)\]
 which by item 5 is \[A\frac{m!}{(i-1)!}(m-i+2)\left(\frac{\lambda}{1-\lambda}\right)^{m-i+1}+O(1/k)\]
On the other hand, by (17), the coefficient of $N^mc$ in $P_{i-1}/k^{i-1}$ can be written as
\[-A\frac{m!}{(i-1)!}(m-i+1)\left(\frac{\lambda}{1-\lambda}\right)^{m-i+1}+O(1/k)\]
Therefore the coefficient of $N^mc, i\le m\le t-1,$ in $\frac{1}{k^{i-1}}(P_{i-1}+\cdots +P_{t-1})$ is
\[A\frac{m!}{(i-1)!}\left(\frac{\lambda}{1-\lambda}\right)^{m-i+1}+O(1/k),\] which approaches to \[\frac{\lambda^{-i+1}}{(i-1)!(1-\lambda)^{m-i+3}}\] as $k\to \infty$ since $A=A(k,m)\to \lambda^{-m}/(m!(1-\lambda)^2)$ as $k\to\infty$. Since $D(k,i-1,\lambda)\to 1$ , $\frac{C(k,i-1)}{k^{i-1}}\to \frac{1}{(i-1)!}$, and since $N^{i-1}$ only appears in $P_{i-1}$, we see that the above sentence is also valid for $m=i-1$.  It then follows from (17) and (18) that the coefficient of $k^{i-1}$ in $H$ is
\[\frac{\lambda^{-i+1}}{(i-1)!}\left[\sum_{j=i-1}^{t-1}\frac{1}{(1-\lambda)^{j-i+3}}N^jc-\frac{1}{1-\lambda}N^{i-1}x\right]\]
Therefore $H$ is of degree $i-2$ or less if and only if
\[N^{i-1}x=\sum_{j=i-1}^{t-1}\frac{1}{(1-\lambda)^{j-i+2}}N^jc\]
The last part of the theorem corresponds to $i=1$ in the above equation. This completes the proof of Theorem 7.

\end{enumerate}
{\bf Appendix}\\ \\
We prove by induction that
\[T_j=
\frac{k}{(1-\lambda)^{j+1}}-\frac{j+1}{(1-\lambda)^{j+2}}+\frac{C(k,j)\lambda^{k-j}}{(1-\lambda)^2}D(k,j,\lambda)\] for $j=0,\cdots,t-1$, where \[D(k,j,\lambda)=\frac{1}{(1-\lambda)^j}[B_0(k,j)\lambda^j+B_1(k,j)\lambda^{j-1}+\cdots+B_{j-1}(k,j)\lambda+1]\] and \[
B_i(k,j)=(-1)^{j-i}C(j,i)\frac{C(k-j,2)}{C(k-i,2)},\
i=0,\cdots,j\]
As we noted in the paper proper this is true for $j=0$ since $D(k,0,\lambda)=1$. By item 1 below $T_j$ can be rewritten as
\[\frac{k}{(1-\lambda)^{j+1}}-\frac{j+1}{(1-\lambda)^{j+2}}+\frac{\lambda^{k-j}A}{(1-\lambda)^{j+2}}\] where
\[A=\sum_{i=0}^j(-1)^{j-i}C(k,i)C(k-i-2,j-i)\lambda^{j-i}\]
Taking the derivative of $T_j$ with respect to $\lambda$, we find
\[\frac{d}{d\lambda}T_j=\frac{k(j+1)}{(1-\lambda)^{j+2}}-\frac{(j+1)(j+2)}{(1-\lambda)^{j+3}}+\frac{\lambda^{k-j-1}([
(1-\lambda)(k-j)+(j+2)\lambda]A+(1-\lambda)\lambda A')}{(1-\lambda)^{j+3}}\]
Then item 3 below proves that
\[\frac{d}{d\lambda}T_j=\frac{k(j+1)}{(1-\lambda)^{j+2}}-\frac{(j+1)(j+2)}{(1-\lambda)^{j+3}}+(j+1)\frac{C(k,j+1)\lambda^{k-j-1}}{(1-\lambda)^2}D(k,j+1,\lambda)\]
This completes the induction since, as stated in the paper proper,
\[T_{j+1}=\frac{1}{j+1}\frac{d}{d\lambda}T_j\]
\begin{enumerate}
\item \[C(k,j)C(j,i)\frac{C(k-j,2)}{C(k-i,2)}=C(k,i)C(k-i-2,j-i)\]
{\bf Proof.}\\
\begin{eqnarray*}& &C(k,j)C(j,i)\frac{C(k-j,2)}{C(k-i,2)}\\
&=&\frac{k!j!(k-j)!(k-i-2)!2!}{j!(k-j)!(j-i)!i!(k-j-2)!2!(k-i)!}\\
&=&\frac{k!(k-i-2)!}{i!(k-i)!(k-j-2)!(j-i)!}\\
&=&C(k,i)C(k-i-2,j-i)
\end{eqnarray*}
Q.E.D.
\item For $k\ge i+2,j\ge i-1$,
\begin{eqnarray*}& &(k-i+1)C(k,i-1)C(k-i-1,j-i+1)+(k-i-j-2)C(k,i)C(k-i-2,j-i)\\&=&(j+1)C(k,i)C(k-i-2,j-i+1).\end{eqnarray*}
{\bf Proof.}\\
\begin{eqnarray*}& &(k-i+1)C(k,i-1)C(k-i-1,j-i+1)+(k-i-j-2)C(k,i)C(k-i-2,j-i)\\
&=&(k-i+1)\frac{k!(k-i-1)!}{(k-i+1)!(i-1)!(k-j-2)!(j-i+1)!}\\& &+(k-i-j-2)\frac{k!(k-i-2)!}{(k-i)!i!(k-j-2)!(j-i)!}\\
&=&\frac{ik!(k-i-1)!}{(k-i)!i!(k-j-2)!(j-i+1)!}\\& &+(k-i-j-2)\frac{(j-i+1)k!(k-i-2)!}{(k-i)!i!(k-j-2)!(j-i+1)!}\\
&=&\frac{k!(k-i-2)![i(k-i-1)+(k-i-j-2)(j-i+1)]}{(k-i)!i!(k-j-2)!(j-i+1)!}\\
&=&\frac{k!(k-i-2)![kj+k-j^2-3j-2]}{(k-i)!i!(k-j-2)!(j-i+1)!}\\
&=&\frac{k!(k-i-2)!(j+1)(k-j-2)}{(k-i)!i!(k-j-2)!(j-i+1)!}\\
&=&\frac{k!(k-i-2)!(j+1)}{(k-i)!i!(k-j-3)!(j-i+1)!}\\
&=&(j+1)C(k,i)C(k-i-2,j-i+1)
\end{eqnarray*}
Q.E.D.

\item Let \[A=\sum_{i=0}^j(-1)^{j-i}C(k,i)C(k-i-2,j-i)\lambda^{j-i}\] and let $A'=\frac{d}{d\lambda}A$. Then
\begin{eqnarray}& &[
(1-\lambda)(k-j)+(j+2)\lambda]A+(1-\lambda)\lambda A'\\
&=& (j+1)\sum_{i=0}^{j+1}(-1)^{j+1-i}C(k,i)C(k-i-2,j+1-i)\lambda^{j+1-i}
\end{eqnarray}
{\bf Proof.} We have \[A'=\sum_{i=0}^{j-1}(-1)^{j-i}(j-i)C(k,i)C(k-i-2,j-i)\lambda^{j-i-1}\]
It is easy to see that
\begin{eqnarray*}& & (1-\lambda)[(k-j)A+\lambda A']\\
&=&(1-\lambda)\sum_{i=0}^j(-1)^{j-i}(k-i)C(k,i)C(k-i-2,j-i)\lambda^{j-i}
\end{eqnarray*}
and
\begin{eqnarray*}
& &-\lambda\sum_{i=0}^j(-1)^{j-i}(k-i)C(k,i)C(k-i-2,j-i)\lambda^{j-i}+(j+2)\lambda A\\
&=&\sum_{i=0}^j (-1)^{j-i+1}(k-i-j-2)C(k,i)C(k-i-2,j-i)\lambda^{j-i+1}.
\end{eqnarray*}
Thus (19) is the sum of \[B=\sum_{i=0}^j(-1)^{j-i}(k-i)C(k,i)C(k-i-2,j-i)\lambda^{j-i}\] and \[C=\sum_{i=0}^j (-1)^{j-i+1}(k-i-j-2)C(k,i)C(k-i-2,j-i)\lambda^{j-i+1}.\]
The coefficient of $\lambda^{j+1-i}$ in $B+C$ is \[(-1)^{j-i+1}[(k-i+1)C(k,i-1)C(k-i-1,j-i+1)+(k-i-j-2)C(k,i)C(k-i-2,j-i)],\]
valid even when $i=0$ or $j+1$ since by convention $C(n,x)=0$ for $x<0$. This is equal to the coefficient of $\lambda^{j+1-i}$ in (20) by item 2 above. Q.E.D.

\end{enumerate}
Acknowledgement: The author wishes to thank Professors Sing-Cheong Ong and Timothy Murphy for their discussions and Professors Roger Horn and Charles Johnson for their communications.


\begin{thebibliography}{99}
\bibitem{1.}  Dowson H. R., Spectral Theory of Linear Operators, Academic Press, 1978.
\bibitem{2.} Lim T.C., Nonexpansive Matrices with Applications to Solutions of Linear Systems
by Fixed Point Iterations, Fixed Point Theory and Applications.Volume 2010 (2010),
Article ID 821928, 13 pages. doi:10.1155/2010/821928.
http://www.hindawi.com/journals/fpta/2010/821928.html
\bibitem{3.} Merris R., Combinatorics, Wiley-Interscience, 2003.
\bibitem{4.} Bayart F. and Matheron E., Dynamics of Linear Operators, Cambridge Tracts in Math., no. 179, 2009.

\end{thebibliography}
\end{document}